\input amstex
\input amsppt.sty
\magnification=\magstep1
\hsize=30truecc
\vsize=22.2truecm
\baselineskip=16truept
\NoBlackBoxes
\TagsOnRight \pageno=1 \nologo
\def\Z{\Bbb Z}
\def\N{\Bbb N}

\def\Q{\Bbb Q}

\def\l{\left}
\def\r{\right}
\def\bg{\bigg}
\def\({\bg(}
\def\[{\bg\lfloor}
\def\){\bg)}
\def\]{\bg\rfloor}
\def\t{\text}
\def\f{\frac}

\def\em{\emptyset}

\def\bi{\binom}
\def\eq{\equiv}

\def\ls{\leqslant}
\def\gs{\geqslant}
\def\mo{\roman{mod}}

\def\da{\delta}

\def\Proof{\noindent{\it Proof}}

\def\Remark{\medskip\noindent{\it  Remark}}

\def\Ack{\medskip\noindent {\bf Acknowledgments}}
\hbox {Preprint, {\tt arXiv:1005.1054}.}
\bigskip
\topmatter
\title On divisibility concerning binomial coefficients\endtitle
\author Zhi-Wei Sun\endauthor
\leftheadtext{Zhi-Wei Sun}
 \rightheadtext{On divisibility concerning binomial coefficients}
\affil Department of Mathematics, Nanjing University\\
 Nanjing 210093, People's Republic of China
  \\  zwsun\@nju.edu.cn
  \\ {\tt http://math.nju.edu.cn/$\sim$zwsun}
\endaffil
\abstract Let $k,l$ and $n$ be positive integers. We mainly show that
$$\align (ln+1)\ \big|&\ k\bi{kn+ln}{kn},
\\2\bi{kn}n\ \big|&\ \bi{2n}{n}C_{2n}^{(k-1)},
\\ \bi{kn}n\ \big|&\ (2k-1)C_n\bi{2kn}{2n},
\\\bi{2n}n\ \big|&\ (k+1)C_n^{(k-1)}\bi{2kn}{kn},
\\2^{k-1}\bi{2n}n\ \big|&\ \bi{2(2^k-1)n}{(2^k-1)n}C_n^{(2^k-2)},
\\(6n+1)\bi{5n}n\ \big|&\ \bi{3n-1}{n-1}C_{3n}^{(4)},
\\\bi{3n}n\ \big|&\ \bi{5n-1}{n-1}C_{5n}^{(2)},
\endalign$$ where $C_n$ denotes
the usual Catalan number $\f1{n+1}\bi{2n}n=\bi{2n}n-\bi{2n}{n+1}$,
and $C_m^{(h)}$ refers to the Catalan number
$\bi{(h+1)m}m/(hm+1)$ of order $h$.
\endabstract
\thanks 2010 {\it Mathematics Subject Classification}.\,Primary 11B65;
Secondary 05A10,  11A07.
\newline\indent Supported by the National Natural Science
Foundation (grant 10871087) and the Overseas Cooperation Fund (grant 10928101) of China.
\endthanks
\endtopmatter
\document

\heading{1. Introduction}\endheading

There are many curious congruences on sums of binomial coefficients (see [S09, S10a, S10b] and [ST1, ST2]).
In a recent paper [S10c] the author investigated products and sums divisible by central binomial coefficients.
(See also [C], [CP] and [GJZ] for another family of sums divisible by central binomial coefficients.)
In this paper we focus on factors of products of two binomial coefficients.

In 2009 J. W. Bober [B] determined all those
$a_1,\ldots,a_r,b_1,\ldots,b_{r+1}\in\Z^+=\{1,2,3,\ldots\}$
with $a_1+\cdots+a_r=b_1+\cdots+b_{r+1}$ such that
$$\f{(a_1n)!\cdots (a_rn)!}{(b_1n)!\cdots(b_{r+1}n)!}$$
is an integer for any $n\in\Z^+$. In particular,
if $k$ and $l$ are positive integers then
$$\align &\f{\bi{ln}n\bi{kln}{ln}}{\bi{kn}n}=\f{(kln)!((k-1)n)!}{(kn)!((l-1)n)!((k-1)ln)!}\in\Z\ \ \t{for all}\ n\in\Z^+,
\\&\ \ \iff k=l,\  \t{or}\ \{k,l\}\cap\{1,2\}\not=\em,\ \t{or}\ \{k,l\}=\{3,5\}.
\endalign$$

Before stating our theorems we need to introduce  Catalan numbers and extended Catalan numbers.

For $n\in\N=\{0,1,2,\ldots\}$, the $n$th (usual) Catalan number is given by
$$C_n=\f1{n+1}\bi{2n}n=\bi{2n}n-\bi{2n}{n-1},$$
The Catalan numbers arise naturally in many enumeration problems in discrete mathematics
(see, e.g., [St, pp.\,219--229]); for example, $C_n$ is the number of
binary parenthesizations of a string of $n+1$ letters.

For $h,n\in\N$ the $n$th (generalized) Catalan number of order $h$ is defined by
$$C_n^{(h)}=\f1{hn+1}\bi{(h+1)n}n=\bi{(h+1)n}n-h\bi{(h+1)n}{n-1}.$$
Note that
$$nC_n^{(h)}=\bi{(h+1)n}{n-1}.$$

Now we present our main results.

\proclaim{Theorem 1.1} {\rm (i)} Let $m\in\N$ and $n\in\Z^+=\{1,2,3,\ldots\}$. Then
$$2\bi{m+n}n\ \bigg|\ \bi{2n}n\bi{2m+2n}{2n},\tag1.1$$
and $\bi{2n}n\bi{2m+2n}{2n}/(2\bi{m+n}n)$ is odd if and only if $n$ is a power of two.
When $n>1$, we have
$$8\bi{m+n}n\ \bigg|\ \bi{2n}n\bi{2m+2n}{2n-1},\tag1.2$$
and
$\bi{2n}n\bi{2m+2n}{2n-1}/(8\bi{m+n}n)$ is odd if and only if $n-1$ is a power of two.

{\rm (ii)} Let $k,n\in\Z^+$. Then
$$2\bi{kn}n\ \bigg|\ \bi{2n}nC_{2n}^{(k-1)}.\tag1.3$$
Moreover, $\bi{2n}nC_{2n}^{(k-1)}/(2\bi{kn}n)$ is odd if and only if $n$ is a power of two.
\endproclaim

\Remark\ 1.1.  After reading [S10c], on May 5, 2010 I. Gessel informed the author that for any $m,n\in\N$ the number
$\bi{2m+2n}{m+n}\bi{m+n}n/\bi{2n}n$ is a positive integer which has a combinatorial interpretation.
Having seen (1.1) in a previous version of this paper, on May 19, 2010 T. Amdeberhan pointed out that
$$\f{\bi{2n}n\bi{2m+2n}{2n}}{\bi{m+n}n}=\f{\bi{2m+2n}{m+n}\bi{m+n}n}{\bi{2m}m}\quad\t{for any}\ m,n\in\N.$$
We note that
$$\f{\bi{2n}n\bi{2m+2n}{2n-1}}{\bi{m+n}{n}}=\f{2\bi{2m+2n}{m+n}\bi{m+n}{n-1}}{\bi{2m+1}m}
\quad\t{for all}\ m\in\N\ \t{and}\ n\in\Z^+.$$

\proclaim{Theorem 1.2} Let $k,n\in\Z^+$.

{\rm (i)} We have
$$\bi{kn}n\ \bigg|\ (2k-1)C_n\bi{2kn}{2n},\tag1.4$$
and $(2k-1)C_n\bi{2kn}{2n}/\bi{kn}n$ is odd if and only if $n+1$ is a power of two.
\smallskip

{\rm (ii)} Let $(k+1)'$ be the odd part of $k+1$. Then
$$\bi{2n}n\ \bigg|\ (k+1)'C_n^{(k-1)}\bi{2kn}{kn},\tag1.5$$
and $(k+1)'C_n^{(k-1)}\bi{2kn}{kn}/\bi{2n}n$ is odd if and only if $(k-1)n+1$ is a power of two.

{\rm (iii)} We also have
$$2^{k-1}\bi{2n}n\ \bigg|\ \bi{2(2^k-1)n}{(2^{k}-1)n}C_n^{(2^{k}-2)}.\tag1.6$$
\endproclaim
\Remark\ 1.2. (1.6) in the case $k=2$ was proved in [S10c].
For a given $k\in\Z^+$
it is interesting to investigate arithmetical properties of the integer sequence
$$S_n^{(k)}:=\f{\bi{2(2^k-1)n}{(2^{k}-1)n}C_n^{(2^{k}-2)}}{2^{k-1}\bi{2n}n}=\f{\bi{2(2^k-1)n}{(2^{k}-1)n}
\bi{(2^k-1)n}n}{2^{k-1}((2^k-2)n+1)\bi{2n}n}\ \ (n=1,2,3,\ldots).$$
\medskip

A key step in our proof of (1.6) is to show the first assertion in our
 following conjecture with $m$ a prime.

\proclaim{Conjecture 1.1} Let $m>1$ be an integer and $k$ and $n$ be positive integers. Then
the sum of all digits in the expansion of $(m^k-1)n$ in base $m$
is at least $k(m-1)$. Also, the expansion of $\f{m^k-1}{m-1}n$ in base $m$
has at least $k$ nonzero digits.
\endproclaim

\proclaim{Theorem 1.3} For any $n\in\Z^+$ we have
$$(6n+1)\bi{5n}n\ \bigg|\ \bi{3n-1}{n-1} C_{3n}^{(4)}\tag1.7$$
and
$$ \bi{3n}n\ \bigg|\ \bi{5n-1}{n-1} C_{5n}^{(2)}.\tag1.8$$
\endproclaim

Define two new sequences $\{s_n\}_{n\gs1}$ and $\{t_n\}_{n\gs1}$ of integers by
$$s_n=\f{\bi{3n-1}{n-1}C_{3n}^{(4)}}{(6n+1)\bi{5n}n}=\f{\bi{3n-1}{n-1}\bi{15n}{3n}}{(6n+1)(12n+1)\bi{5n}n}\tag1.9$$
and
$$t_n=\f{\bi{5n-1}{n-1}C_{5n}^{(2)}}{\bi{3n}n}=\f{\bi{5n-1}{n-1}\bi{15n}{5n}}{(10n+1)\bi{3n}n}.\tag1.10$$
Then the values of $s_1,\ldots,s_8$ are
$$\gather 1,\ 203,\ 77572,\ 38903007,\ 22716425576,
\\14621862696188,\ 10071456400611060,\ 7291908546474763815
\endgather$$
respectively, while the values of $t_1,\ldots,t_5$ are
$$ 91,\ 858429,\ 12051818636,\ 200142760587609,\ 3648677478873075576$$
respectively.
It would be interesting to find recursion formulae or combinatorial interpretations for $s_n$ and $t_n$.

Based on our computation via {\tt Mathematica}, we formulate the following conjecture on the sequence $\{t_n\}_{n\gs1}$.

\proclaim{Conjecture 1.2} Let $n$ be any positive integer. We have
$$21t_n\eq0\ (\mo\ 10n+3).\tag1.11$$
If $3\nmid n$, then $(10n+3)\mid7t_n$.
If $7\nmid n+1$, then $(10n+3)\mid 3t_n$.
\endproclaim

In view of Theorems 1.1-1.3 and some computation,  we raise the following conjecture.

\proclaim{Conjecture 1.3} Let $k$ and $l$ be integers greater than one.
If $\bi{kn}n\mid\bi{ln}n\bi{kln}{ln-1}$ for all $n\in\Z^+$, then $k=l$, or $l=2$, or $\{k,l\}=\{3,5\}$.
If $\bi{kn}n\mid\bi{ln}{n-1}\bi{kln}{ln}$ for all $n\in\Z^+$, then $k=2$, and $l+1$ is a power of two.
\endproclaim

Recall that $\bi{hn+n}n/(hn+1)=C_n^{(h)}\in\Z$ for all $h,n\in\Z^+$.
It is natural to ask for what positive integers $k$ and $l$ we have
$(ln+1)\mid\bi{kn+ln}{kn}$ for all $n\in\Z^+$. In this direction, we obtain the following result.

\proclaim{Theorem 1.4} Let $k,l,n\in\Z^+$. Then
$$\bi{kn+ln}{kn}\eq0\ \l(\mo\ \f{ln+1}{(k,ln+1)}\r),\tag1.12$$
where $(k,ln+1)$ denotes the greatest common divisor of $k$ and $ln+1$.
In particular, $(ln+1)\mid\bi{kn+ln}{kn}$ if $l$ is divisible by all prime factors of $k$.
\endproclaim

Our following conjecture seems difficult.

\proclaim{Conjecture 1.4} Let $k$ and $l$ be positive integers.
If $(ln+1)\mid\bi{kn+ln}{kn}$ for all sufficiently large positive integers $n$, then
each prime factor of $k$ divides $l$. In other words, if $k$ has a prime factor not dividing $l$ then
there are infinitely many positive integers $n$ such that  $(ln+1)\nmid\bi{kn+ln}{kn}$.
\endproclaim

In view of Conjecture 1.4 we introduce a new function $f:\Z^+\times\Z^+\to \N$.
For positive integers $k$ and $l$, if $(ln+1)\mid\bi{kn+ln}{kn}$ for all $n\in\Z^+$
(which happens if all prime factors of $k$ divide $l$) then we set $f(k,l)=0$, otherwise
we define $f(k,l)$ to be the smallest positive integer $n$ such that $(ln+1)\nmid\bi{kn+ln}{kn}$.
Our computation via {\tt Mathematica} yields the following values of $f$:

$$\gather f(7,36)=279,\ f(10,192)=362,\ f(11,100)=1187,
\\ f(13,144)=2001,\ f(22,200)=6462,\ f(31,171)=1765;
\\f(43,26)=640,\ f(53,32)=790,\ f(67,56)=2004,
\\f(73,61)=2184,\ f(74,62)=885,\ f(97,81)=2904.
\endgather$$
It would be interesting to investigate the behavior of the function $f$.

In the next section we will establish three auxiliary theorems on inequalities
involving the floor function. Sections 3 is devoted to the proofs of
Theorems 1.1--1.4.

Throughout this paper, for a real number $x$ we let $\{x\}=x-\lfloor
x\rfloor$ be the fractional part of $x$.

\heading{2. Three auxiliary theorems on inequalities}\endheading

\proclaim{Theorem 2.1}  Let $m\in\Z^+$ and
$k,n\in\Z$. Then we have
$$\aligned&\l\lfloor\f{2kn}m\r\rfloor-\l\lfloor\f{kn}m\r\rfloor
+\l\lfloor\f{(k-1)n}m\r\rfloor-\l\lfloor\f{2(k-1)n}m\r\rfloor
\\&\ \ \gs\l\lfloor\f{n+1}m\r\rfloor-\l\lfloor\f{2k-1}m\r\rfloor+\l\lfloor\f{2k-2}m\r\rfloor,
\endaligned\tag2.1$$
unless $2\mid m$, $k\eq m/2+1\ (\mo\ m)$ and $n\eq-1\ (\mo\ m)$,
in which case the left-hand side of (2.1) minus the right-hand side
equals $-1$.
\endproclaim
\Proof. Clearly (2.1) holds when $m=1$. Below we assume that $m\gs2$.

Let $A_m(k,n)$ denote the left-hand side of (2.1) minus the
right-hand side. Then
$$\align -A_m(k,n)=&\l\{\f{2kn}m\r\}-\l\{\f{kn}m\r\}+\l\{\f{(k-1)n}m\r\}-\l\{\f{2(k-1)n}m\r\}
\\&-\l\{\f {n+1}m\r\}+\l\{\f{2k-1}m\r\}-\l\{\f{2k-2}m\r\}.
\endalign$$
Hence $A_m(k,n)\gs0$ if and only if
$$\l\{\f{2kn}m\r\}-\l\{\f{kn}m\r\}+\l\{\f{(k-1)n}m\r\}-\l\{\f{2(k-1)n}m\r\}+\l\{\f{2k-1}m\r\}-\l\{\f{2k-2}m\r\}<1.\tag2.2$$
(Note that $2kn-kn+(k-1)n-2(k-1)n+(2k-1)-(2k-2)=n+1$.)

{\it Case} 1. $\{kn/m\}<1/2\ \&\ \{(k-1)n/m\}<1/2$, or ($\{kn/m\}\gs1/2\ \&\ \{(k-1)n/m\}\gs1/2$).

In this case, the left-hand side of (2.2) equals
$$C:=\l\{\f{kn}m\r\}-\l\{\f{(k-1)n}m\r\}+\l\{\f{2k-1}m\r\}-\l\{\f{2k-2}m\r\}.$$
If $m\nmid (k-1)n$, then
$$C<\l\{\f{kn}m\r\}+\f1m\ls 1.$$
If $m\mid(k-1)n$ and $n\not\eq-1\ (\mo\ m)$, then
$$C\ls \l\{\f nm\r\}+\f 1m<1.$$
If $m\mid(k-1)n$ and $n\eq-1\ (\mo\ m)$, then $\{kn/m\}=(m-1)/m\gs1/2>\{(k-1)n/m\}=0$
which leads a contradiction.

{\it Case} 2. $\{kn/m\}<1/2\ls\{(k-1)n/m\}$.

 In this case, the left-hand side of (2.2) equals
 $$D:=\l\{\f{kn}m\r\}-\l\{\f{(k-1)n}m\r\}+1+\l\{\f{2k-1}m\r\}-\l\{\f{2k-2}m\r\}.$$
If $n\not\eq-1\ (\mo\ m)$, then $\{(k-1)n/m\}-\{kn/m\}\not=1/m$ and hence
 $$D<-\f1m+1+\f1m=1.$$
 If $n\eq-1\ (\mo\ m)$ and $2k\eq 1\ (\mo\ m)$, then
 $$D=-\f1m+1+0-\f{m-1}m<1.$$
 If $n\eq-1\ (\mo\ m)$ and $2k\not\eq1\ (\mo\ m)$, then we must have $2\mid m$ and $k\eq m/2+1\ (\mo\ m)$
 since
 $$\l\{\f{-k}m\r\}<\f12\ls\l\{\f{1-k}m\r\}.$$

When $2\mid m$, $k\eq m/2+1\ (\mo\ m)$ and $n\eq-1\ (\mo\ m)$, it is easy to see that
$$-A_m(k,n)=\f{m-2}m-\f{m/2-1}m+\f{m/2}m+\f1m=1.$$

 {\it Case} 3. $\{kn/m\}\gs1/2>\{(k-1)n/m\}$.

In this case, the left-hand side of (2.2) is
$$\l\{\f{kn}m\r\}-1-\l\{\f{(k-1)n}m\r\}+\l\{\f{2k-1}m\r\}-\l\{\f{2k-2}m\r\}\ls\l\{\f{kn}m\r\}-1+\f1m\ls0.$$

Combining the above we have completed the proof of Theorem 2.1. \qed

\proclaim{Theorem 2.2}  Let $m>2$ be an integer. For any $k,n\in\Z$ we have
$$\l\lfloor\f{2kn}m\r\rfloor+\l\lfloor\f{n}m\r\rfloor+\l\lfloor\f{k+1}m\r\rfloor
\gs\l\lfloor\f km\r\rfloor+\l\lfloor\f{2n}m\r\rfloor+\l\lfloor\f{kn}m\r\rfloor+\l\lfloor\f{(k-1)n+1}m\r\rfloor.\tag2.3$$
\endproclaim
\Proof. As $\lfloor x\rfloor=x-\{x\}$ for all $x\in\Q$, (2.3) is equivalent to the inequality $M\gs\{(k+1)/m\}$, where
$$M:=\l\{\f km\r\}+\l\{\f{(k-1)n+1}m\r\}+\l\{\f{kn}m\r\}-\l\{\f{2kn}m\r\}+\l\{\f{2n}m\r\}-\l\{\f nm\r\}.$$
Since $k+((k-1)n+1)+kn-2kn+2n-n=k+1$, it suffices to show that $M\gs0$.

{\it Case} 1. $\{n/m\}<1/2\ \&\ \{kn/m\}<1/2$, or $(\{n/m\}\gs1/2\ \&\ \{kn/m\}\gs1/2)$.

 In this case,
$$M=\l\{\f km\r\}+\l\{\f{(k-1)n+1}m\r\}+\l\{\f nm\r\}-\l\{\f{kn}m\r\}.$$
If $kn\eq-1\ (\mo\ m)$, then $m\nmid n$ and hence
$$M=\l\{\f km\r\}+\l\{\f{-n}m\r\}+\l\{\f nm\r\}-\l\{\f{-1}m\r\}=\l\{\f km\r\}+1-\f{m-1}m>0.$$
If $kn\not\eq-1\ (\mo\ m)$, then
$$M\gs\l\{\f km\r\}+\l\{\f{kn+1}m\r\}-\l\{\f{kn}m\r\}=\l\{\f km\r\}+\f 1m>0.$$

 {\it Case} 2. $\{n/m\}<1/2\ls\{kn/m\}$.

In this case,
$$M=\l\{\f km\r\}+\l\{\f{(k-1)n+1}m\r\}+1-\l\{\f{kn}m\r\}+\l\{\f nm\r\}>0.$$

{\it Case} 3. $\{kn/m\}<1/2\ls\{n/m\}$.

In this case, $m\nmid n$ and
$$M=\l\{\f km\r\}+\l\{\f{(k-1)n+1}m\r\}-\l\{\f{kn}m\r\}+\l\{\f nm\r\}-1.$$

If $(k-1)n+1\eq0\ (\mo\ m)$, then $\{(n-1)/m\}=\{kn/m\}<1/2\ls\{n/m\}$, hence
$m$ is odd (otherwise $n\eq m/2\ (\mo\ m)$ and thus $1\eq0\ (\mo\ m/2)$ which is impossible)
and $n\eq (m+1)/2\ (\mo\ m)$, therefore $k-1\eq (k-1)2n\eq-2\ (\mo\ m)$ and
$$M=\l\{\f km\r\}-\l\{\f{n-1}m\r\}+\l\{\f nm\r\}-1=\l\{\f km\r\}-\f{m-1}m=0.$$

Now assume that $(k-1)n+1\not \eq0\ (\mo\ m)$. Clearly $\{kn/m\}<\{(n-1)/m\}$ and hence
$$M=\l\{\f km\r\}+\l(\l\{\f{kn}m\r\}-\l\{\f{n-1}m\r\}+1\r)-\l\{\f{kn}m\r\}+\l\{\f nm\r\}-1
\gs\f 1m.$$

By the above we always have $M\gs0$ and hence (2.3) follows. \qed

\proclaim{Lemma 2.1} {\rm (i)} For any real number $x$ we have
$$\{12x\}+\{5x\}+\{2x\}\gs\{4x\}+\{15x\}.\tag2.4$$

{\rm (ii)} Let $x$ be a real number with $\{5x\}\gs\{2x\}\gs1/2$. Then $\{5x\}\gs2/3$.
\endproclaim
\Proof. (i) Since $12x+5x+2x-4x=15x$, (2.4) reduces to
$$X:=\{12x\}+\{5x\}+\{2x\}-\{4x\}\gs0.$$

If $\{2x\}\gs 1/2$, then $X>0$ since
$$\{2x\}-\{4x\}=\{2x\}-(2\{2x\}-1)=1-\{2x\}>0.$$
If $\{2x\}<1/2$ and $\{2x\}\ls \{5x\}$, then
$$X=\{12x\}+\{5x\}-\{2x\}\gs0.$$

 Below we assume that $\{5x\}<\{2x\}<1/2$. Clearly $\{3x\}=\{5x\}-\{2x\}+1$.

 {\it Case} 1. $\{2x\}=2\{x\}$ and hence $\{x\}<1/4$.

Since $\{5x\}<2\{x\}$, we cannot have $\{5x\}=5\{x\}$. As $5\{x\}<5/4<2$, we must have $\{5x\}=5\{x\}-1$
and hence $\{x\}\gs1/5$. Note that
$12/5\ls12\{x\}<3$ and hence $\{12x\}=12\{x\}-2$. Therefore
$$X=12\{x\}-2+5\{x\}-1-2\{x\}=15\{x\}-3\gs0.$$

{\it Case} 2. $\{2x\}=2\{x\}-1$ and hence $1/2\ls \{x\}<3/4$.

As $5/2\ls 5\{x\}<15/4$, \ $\{5x\}$ is $5\{x\}-2$ or $5\{x\}-3$.
Since $5\{x\}-2>\{2x\}=2\{x\}-1$, we must have $\{5x\}=5\{x\}-3<\{2x\}=2\{x\}-1$ and hence $3/5\ls\{x\}<2/3$.
Observe that
$$7<12\times \f 35\ls 12\{x\}<12\times\f 23=8$$
and thus $\{12x\}=12\{x\}-7$. Therefore
$$X=12\{x\}-7+5\{x\}-3-(2\{x\}-1)=15\{x\}-9\gs0.$$

Combining the above we have proved (2.4).

\medskip

(ii) Since $\{3x\}=\{5x\}-\{2x\}<1-\{2x\}\ls1/2\ls\{2x\}$, we have $\{3x\}=\{2x\}+\{x\}-1$
and hence
$$\{5x\}=\{3x\}+\{2x\}=2\{2x\}+\{x\}-1.$$

{\it Case} 1. $\{2x\}=2\{x\}$ and hence $1/4\ls \{x\}<1/2$.

In this case, $\{3x\}=3\{x\}-1$ and
$$\{5x\}=5\{x\}-1\gs\f 53-1=\f 23.$$

{\it Case} 2. $\{2x\}=2\{x\}-1$ and hence $\{x\}\gs3/4$.

In this case,
$$\{5x\}=2(2\{x\}-1)+\{x\}-1=5\{x\}-3\gs\f{15}4-3=\f34>\f23.$$

So far we have also completed the proof of the second part of Lemma 2.1. \qed

\proclaim{Theorem 2.3} Let $m>1$ and $n$ be integers.

{\rm (i)} If $3\nmid m$, then
$$\l\lfloor\f{15n-1}m\r\rfloor+\l\lfloor\f{2}m\r\rfloor+\l\lfloor\f{4n}m\r\rfloor
\gs\l\lfloor\f{12n+2}m\r\rfloor+\l\lfloor\f{2n}m\r\rfloor+\l\lfloor\f{5n-1}m\r\rfloor.\tag2.5$$

{\rm (ii)} If $5\nmid m$, then
$$\l\lfloor\f{15n-1}m\r\rfloor+\l\lfloor\f{2n}m\r\rfloor
\gs\l\lfloor\f{10n+1}m\r\rfloor+\l\lfloor\f{4n}m\r\rfloor+\l\lfloor\f{3n-1}m\r\rfloor.\tag2.6$$
\endproclaim

\Remark. For a positive integer $m$ divisible by 3, we can prove that (2.5) holds unless $n\eq 2m/3\ (\mo\ m)$
in which case the left-hand side of (2.5) minus the right-hand side equals $-1$.
For a positive integer $m$ divisible by 5, we can show that (2.6) holds unless $n\eq 2m/5,4m/5\ (\mo\ m)$
in which case the left-hand side of (2.6) minus the right-hand side equals $-1$.

\medskip
\noindent{\it Proof of Theorem 2.3}.
 (i) Clearly (2.5) holds when $m=2$. Below we assume that $m>2$ and $3\nmid m$.

As $\lfloor x\rfloor=x-\{x\}$ for all $x\in\Q$, (2.5) has the following equivalent form:
$$\l\{\f{12n+2}m\r\}+\l\{\f{5n-1}m\r\}+\l\{\f{2n}m\r\}-\l\{\f{4n}m\r\}\gs\l\{\f{15n-1}m\r\}+\f2m.\tag2.7$$
Since $m\mid 15n$ if and only if $m\mid 5n$, we have
$$\l\{\f{5n-1}m\r\}-\l\{\f{15n-1}m\r\}=\l\{\f{5n}m\r\}-\l\{\f{15n}m\r\}$$
and thus (2.7) can be written as
$$\l\{\f{12n+2}m\r\}+\l\{\f{5n}m\r\}+\l\{\f{2n}m\r\}-\l\{\f{4n}m\r\}\gs\l\{\f{15n}m\r\}+\f2m.\tag2.8$$

{\it Case} 1. $12n+\da\eq0\ (\mo\ m)$ for some $\da\in\{1,2\}$.

In this case, $m$ does not divide $3n$ and (2.7) can be rewritten as
$$\l\{\f{5n}m\r\}+\l\{\f{2n}m\r\}-\l\{\f{4n}m\r\}\gs\l\{\f{3n-\da}m\r\}+\f{\da}m=\l\{\f{3n}m\r\}.$$
(Note that if $m\mid 12n+2$ and $m\mid 3n-1$
then $m$ divides $12n+2-4(3n-1)=6$ which contradicts that $m>2$ and $3\nmid m$.)
Since $5n+2n-4n=3n$, it suffices to prove that
$$\l\{\f{5n}m\r\}+\l\{\f{2n}m\r\}-\l\{\f{4n}m\r\}\gs0.$$

If $\{2n/m\}\gs1/2$, then
$$\l\{\f{2n}m\r\}-\l\{\f{4n}m\r\}=\l\{\f{2n}m\r\}-\(2\l\{\f{2n}m\r\}-1\)=1-\l\{\f{2n}m\r\}>0.$$
So we simply suppose that $\{2n/m\}<1/2$ and want to prove the inequality $\{5n/m\}\gs\{2n/m\}$.

{\it Case} 1.1. $m\eq\da\ (\mo\ 3)$.

In this case, we have $4n\eq(m-\da)/3\ (\mo\ m)$
and hence $2\{2n/m\}=\{4n/m\}<1/3$.
If $\{n/m\}\ls 2/3$, then
$$\l\{\f{5n}m\r\}=\l\{\f{4n}m\r\}+\l\{\f nm\r\}\gs \l\{\f{4n}m\r\}\gs\l\{\f{2n}m\r\}.$$
If $\{n/m\}>2/3$, then $\{n/m\}\gs(2m+\da)/(3m)$ (since $m\eq \da\not\eq0\ (\mo\ m)$), $\{2n/m\}=2\{n/m\}-1>1/3$,
and hence
$$\l\{\f{5n}m\r\}=\l\{\f{4n}m\r\}+\l\{\f nm\r\}-1=2\l\{\f{2n}m\r\}+\l\{\f nm\r\}-1>\l\{\f{2n}m\r\}.$$

{\it Case} 1.2. $m\eq-\da\ (\mo\ 3)$.

In this case, we have $4n\eq-(m+\da)/3\ (\mo\ m)$
and hence $2\{2n/m\}=\{4n/m\}=1-(m+\da)/(3m)=2/3-\da/(3m)$.
If $\{n/m\}\ls 1/3$, then
$$\l\{\f{5n}m\r\}=\l\{\f{4n}m\r\}+\l\{\f nm\r\}\gs \l\{\f{4n}m\r\}\gs\l\{\f{2n}m\r\}.$$
If $\{n/m\}>1/3$, then $\{n/m\}\gs(m+\da)/(3m)$ (since $m\eq -\da\not\eq0\ (\mo\ m)$), $1/2>\{2n/m\}=2\{n/m\}-1$,
hence $3\{n/m\}-2<9/4-2=1/4$ and
$$\l\{\f{5n}m\r\}=\l\{\f{2n}m\r\}+\l\{\f {3n}m\r\}\gs\l\{\f{2n}m\r\}$$
provided $\{3n/m\}=3\{n/m\}-2$.
If $\{n/m\}>1/3$ and $\{3n/m\}\not=3\{n/m\}-2$, then $\{3n/m\}=3\{n/m\}-1$, hence $\{n/m\}<2/3$
and $\{n/m\}\ls(2m-\da)/(3m)=\{4n/m\}$, therefore
$$\l\{\f{5n}m\r\}=\l\{\f{4n}m\r\}+\l\{\f nm\r\}-1\gs\l\{\f{2n}m\r\}=2\l\{\f nm\r\}-1.$$

 Combining our discussions in the cases 1.1 and 1.2, we obtain the desired result in Case 1.

{\it Case}\ 2. $12n+1,12n+2\not\eq0\ (\mo\ m)$.

In this case, (2.8) is equivalent to the inequality
$$\{12x\}+\{5x\}+\{2x\}-\{4x\}\gs\{15x\}$$
with $x=n/m$, which follows from Lemma 2.1(i).

 Combining the above we have proved the first part of Theorem 2.3.

 \medskip

(ii) As $\lfloor x\rfloor=x-\{x\}$ for all $x\in\Q$, (2.6) has the following equivalent form:
$$\l\{\f{10n+1}m\r\}+\l\{\f{3n-1}m\r\}+\l\{\f{4n}m\r\}-\l\{\f{2n}m\r\}\gs\l\{\f{15n-1}m\r\}+\f1m.$$
 Suppose that $5\nmid m$. Then $m\mid 15n$ if and only if $m\mid 3n$. Thus the last inequality can be rewritten as
 $$\l\{\f{10n+1}m\r\}-\f1m+\l\{\f{3n}m\r\}+\l\{\f{4n}m\r\}-\l\{\f{2n}m\r\}\gs\l\{\f{15n}m\r\}$$
 which is equivalent to
$$W:=\l\{\f{10n+1}m\r\}-\f1m+\l\{\f{3n}m\r\}+\l\{\f{4n}m\r\}-\l\{\f{2n}m\r\}\gs0\tag2.9$$
since $(10n+1)-1+3n+4n-2n=15n$.

In the case $m\mid 3n$, (2.9) reduces to
$$\l\{\f{n+1}m\r\}+\l\{\f nm\r\}\gs\l\{\f{2n}m\r\}+\f 1m,$$
which holds since
$$\l\{\f{n+1}m\r\}+\l\{\f nm\r\}\gs\l\{\f{2n+1}m\r\}.$$
(If $m\mid 3n$ and $m\mid 2n+1$, then $m\nmid n$ and hence
$\{(n+1)/m\}+\{n/m\}=\{-n/m\}+\{n/m\}=1=\{2n/m\}+1/m$.)

Below we assume that $m\nmid 3n$. Then
$$W:=\l\{\f{10n+1}m\r\}+\l\{\f{3n-1}m\r\}+\l\{\f{4n}m\r\}-\l\{\f{2n}m\r\}.$$

If $\{2n/m\}<1/2$, then $\{4n/m\}-\{2n/m\}=\{2n/m\}\gs0$.  If $\{2n/m\}\gs1/2$ and $\{(5n-1)/m\}<\{2n/m\}$, then
$$\align W=&\l\{\f{10n+1}m\r\}+\l\{\f{3n-1}m\r\}+\l\{\f{2n}m\r\}-1
\\=&\l\{\f{10n+1}m\r\}+\(\l\{\f{5n-1}m\r\}+1\)-1\gs0.
\endalign$$

Now we consider the remaining case $\{(5n-1)/m\}\gs\{2n/m\}\gs1/2$. Note that
$$W=\l\{\f{10n+1}m\r\}+\l\{\f{3n-1}m\r\}+\l\{\f{2n}m\r\}-1=\l\{\f{10n+1}m\r\}+\l\{\f{5n-1}m\r\}-1.$$
Clearly $W=0$ if $m\mid 5n$. If $m\mid 10n+1$, then $2\nmid m$, $5n\eq(m-1)/2\ (\mo\ m)$ and hence $\{(5n-1)/m\}<1/2$.

Below we simply assume that $m\nmid 5n$ and $m\nmid 10n+1$. Observe that
$$\align W=&\l\{\f{10n}m\r\}+\f1m+\l\{\f{5n-1}m\r\}-1=\l\{\f{10n}m\r\}+\l\{\f{5n}m\r\}-1
\\=&2\l\{\f{5n}m\r\}-1+\l\{\f{5n}m\r\}-1=3\l\{\f{5n}m\r\}-2.
\endalign$$
Set $x=n/m$. Then $\{5x\}\gs\{2x\}\gs1/2$. By Lemma 2.1(ii), $\{5x\}\gs2/3$ and hence $W\gs0$.
This concludes the proof. \qed

\heading{3. Proofs of Theorems 1.1--1.4}\endheading

For a prime $p$,  the $p$-adic evaluation of an integer $m$ is given by
$$\nu_p(m)=\sup\{a\in\N:\ p^a\mid m\}.$$
For a rational number $x=m/n$ with $m\in\Z$ and $n\in\Z^+$, we set $\nu_p(x)=\nu_p(m)-\nu_p(n)$ for any prime $p$.
Note that a rational number $x$ is an integer if and only if $\nu_p(x)\gs0$ for all primes $p$.

Let $p$ be any prime. A useful theorem of Legendre (see, e.g., [R, pp. 22--24]) asserts that for any $n\in\N$ we have
$$\nu_p(n!)=\sum_{i=1}^\infty\l\lfloor\f n{p^i}\r\rfloor=\f{n-\rho_p(n)}{p-1},$$
 where $\rho_p(n)$ is the sum of the digits of $n$ in the expansion of $n$ in base $p$.
This immediately yields the following lemma.

\proclaim{Lemma 3.1} Let $n$ be a positive integer. Then $\nu_2(n!)\ls n-1$. Also, $\nu_2(n!)=n-1$
if and only if $n$ is a power of two.
\endproclaim
\Remark. The latter part of Lemma 3.1 also follows from [SD, Lemma 4.1].

\medskip\noindent
{\it Proof of Theorem 1.1}.  (i) Set
$$Q(m,n):=\f{\bi{2n}n\bi{2m+2n}{2n}}{2\bi{m+n}n}.$$
Then
$$Q(m,n)=\f{\prod_{j=1}^n(2j)(2j-1)}{2(n!)^2}\prod_{j=1}^n\f{2m+2j-1}{2j-1}
=\f{2^{n-1}}{n!}\prod_{j=1}^n(2m+2j-1).$$
If $p$ is an odd prime and $m'$ is an integer with $m'\eq m-1/2\ (\mo\ p^{\nu_p(n!)+1})$, then
$$\f{\prod_{j=1}^n(2m+2j-1)}{n!}\eq\f{2^n\prod_{j=1}^n(m'+j)}{n!}=2^n\bi{m'+n}n\ (\mo\ p).$$
So $Q(m,n)$ is a $p$-adic integer for any odd prime $p$.
Note also that
$$\nu_2(Q(m,n))=n-1-\nu_2(n!).$$
Applying Lemma 3.1 we see that $Q(m,n)\in\Z$ and that $2\nmid Q(m,n)$ if and only if $n$ is a power of two.

When $n>1$ we have
$$\f{\bi{2n}n\bi{2m+2n}{2n-1}}{8\bi{m+n}n}=Q(m+1,n-1).$$
So the latter assertion in Theorem 1.1(i) follows from the above.

(ii) Observe that
$$\f{\bi{2n}nC_{2n}^{(k-1)}}{2\bi{kn}n}=\f{\bi{2n}n}{2\bi{kn}n}\(\bi{2kn}{2n}-(k-1)\bi{2kn}{2n-1}\).$$
So it suffices to apply Theorem 1.1(i) with $m=(k-1)n$.

\medskip

The proof of Theorem 1.1 is now complete. \qed

\proclaim{Lemma 3.2} Let $p$ be a prime and let $k\in\N$ and $n\in\Z^+$. Then
$$\f{\rho_p((p^k-1)n)}{p-1}=\sum_{j=1}^\infty\l\{\f{(p^k-1)n}{p^j}\r\}\gs k\tag3.1$$
and hence the expansion of $(p^k-1)n$ in base $p$ has at least $k$ nonzero digits.
\endproclaim
\Proof. For any $m\in\Z^+$, by Legendre's theorem we have
$$\f{\rho_p(m)}{p-1}=\f m{p-1}-\nu_p(m!)
=\sum_{j=1}^\infty\f m{p^j}-\sum_{j=1}^\infty\l\lfloor\f m{p^j}\r\rfloor=\sum_{j=1}^\infty\l\{\f m{p^j}\r\}.$$
If the expansion of $m$ in base $p$ has less than $k$ nonzero digits, then $\rho_p(m)<k(p-1)$.
So it remains to show the inequality in (3.1).

Observe that
$$p^k\bi{p^kn-1}{n-1}=\bi{p^kn}{n}=\f{(p^kn)!}{n!((p^k-1)n)!}$$
and
$$\align &\nu_p((p^kn)!)-\nu_p(n!)-\nu_p(((p^k-1)n)!)
\\=&\sum_{j=1}^\infty\l\lfloor\f{p^kn}{p^j}\r\rfloor
-\sum_{j=1}^\infty\l\lfloor\f n{p^j}\r\rfloor-\sum_{j=1}^\infty\l\lfloor\f{(p^k-1)n}{p^j}\r\rfloor
\\=&\sum_{j=1}^k p^{k-j}n-\sum_{j=1}^\infty\l\lfloor\f{(p^k-1)n}{p^j}\r\rfloor=\sum_{j=1}^\infty\l\{\f{(p^k-1)n}{p^j}\r\}.
\endalign$$
So the inequality in (3.1) follows. We are done. \qed

\medskip
\noindent{\it Proof of Theorem 1.2}. (i) Define $A_m(k,n)$ for $m>1$ as in the proof of Theorem 2.1.
Observe that
$$Q_1:=\f{(2k-1)C_n\bi{2kn}{kn}}{\bi{kn}n}=\f{(2kn)!((k-1)n)!(2k-1)!}{(n+1)!(kn)!(2(k-1)n)!(2k-2)!}.$$
So, for any prime $p$ we have
$$\nu_p(Q_1)=\sum_{i=1}^\infty A_{p^i}(k,n).$$
By Theorem 2.1, $A_{p^i}(k,n)\gs0$ unless $p=2$, $k\eq 2^{i-1}+1\ (\mo\ 2^i)$ and $n\eq-1\ (\mo\ 2^i)$,
in which case $A_{p^i}(k,n)=-1$.
Therefore $2Q_1\in\Z$.

Observe that
$$\align Q_1=&\f{2k-1}{n+1}\cdot\f{\prod_{j=1}^n(2j)(2j-1)}{(n!)^2}\prod_{j=1}^n\f{(2k-2)n+2j-1}{2j-1}
\\=&\f{2^n(2k-1)}{(n+1)!}\prod_{j=1}^n((2k-2)n+2j-1).
\endalign$$
and thus $\nu_2(Q_1)=n-\nu_2((n+1)!)$. With the help of Lemma 3.1 we obtain that
$Q\in\Z$ and that $Q$ is odd if and only if $n+1$ is a power of two. This proves Theorem 1.2(i).

 (ii) Obviously
 $$Q_2:=\f{(k+1)C_n^{(k-1)}\bi{2kn}{kn}}{\bi{2n}n}
=\f{(k+1)!(2kn)!n!}{k!(kn)!((k-1)n+1)!(2n)!}.$$
Given an odd prime $p$, clearly $\nu_p(Q_2)$ coincides with
$$\sum_{i=1}^\infty\(\l\lfloor\f{k+1}{p^i}\r\rfloor-\l\lfloor\f k{p^i}\r\rfloor
+\l\lfloor\f{2kn}{p^i}\r\rfloor-\l\lfloor\f{kn}{p^i}\r\rfloor+\l\lfloor\f n{p^i}\r\rfloor
-\l\lfloor\f{2n}{p^i}\r\rfloor-\l\lfloor\f{(k-1)n+1}{p^i}\r\rfloor\),$$
which is nonnegative by Theorem 2.2.

 Now we consider $\nu_2(Q_2)$. Set $m=(k-1)n$. Then
 $$\align Q_2=&\f{(k+1)\bi{2m+2n}{m+n}\bi{m+n}n}{(m+1)\bi{2n}n}=\f{(k+1)4^m}{(m+1)!}\prod_{j=1}^m\l(j+n-\f12\r)
 \\=&\f{2^m(k+1)}{(m+1)!}\prod_{j=1}^m(2j+2n-1)
 \endalign$$
and therefore $\nu_2(Q_2)=\nu_2(k+1)+m-\nu_2((m+1)!)$. Applying Lemma 3.1 we see that $\nu_2(Q_2)\gs \nu_2(k+1)$.
So $Q_2/2^{\nu_2(k+1)}$ is an integer. With the help of Lemma 3.1,
$$\align &\f{Q_2}{2^{\nu_2(k+1)}}=\f{(k+1)'C_n^{(k-1)}\bi{2kn}{kn}}{\bi{2n}n}\ \t{is odd}
\\\iff& \ \nu_2((m+1)!)=m
\\\iff&\ m+1=(k-1)n+1\ \t{is a  power of two}.
\endalign$$
This concludes the proof of Theorem 1.2(ii).

  (iii)  Since the odd part of $(2^k-1)+1$ is 1, by part (ii)
$$Q_3:=\f{\bi{2(2^k-1)n}{(2^k-1)n}C_n^{(2^k-2)}}{\bi{2n}n}$$
is an integer and also $\nu_2(Q_3)=m-\nu_2((m+1)!)$, where $m=((2^k-1)-1)n$ is even.
Thus, with helps of Legendre's theorem and Lemma 3.2 with $p=2$, we have
$$\nu_2(Q_3)=m!-\nu_2(m!)=\rho_2(m)=\rho_2((2^{k-1}-1)n)\gs k-1.$$
Therefore $2^{k-1}\mid Q_3$ and hence (1.6) holds.

\medskip

So far we have completed the proof of Theorem 1.2. \qed

\medskip
\noindent{\it Proof of Theorem 1.3}. Observe that
$$A:=\f{\bi{3n-1}{n-1}C_{3n}^{(4)}}{(6n+1)\bi{5n}n}=\f{(15n-1)!2!(4n)!}{(12n+2)!(2n)!(5n-1)!}$$
and
$$B:=\f{\bi{5n-1}{n-1}C_{5n}^{(2)}}{\bi{3n}n}=\f{(15n-1)!(2n)!}{(10n+1)!(4n)!(3n-1)!}.$$
By Theorem 2.3, $\nu_p(A)\gs0$ for any prime $p\not=3$, and
$\mu_p(B)\gs0$ for any prime $p\not=5$. Thus, it suffices to show that
$\nu_3(A)\gs0$ and $\nu_5(B)\gs0$. In fact,
$$\f{C_{3n}^{(4)}}{(6n+1)\bi{5n}n}=\f1{(6n+1)(12n+1)}\prod\Sb j=1\\3\nmid j\endSb^{3n}\f{12n+j}j$$
is a $3$-adic integer, and
$$\f{C_{5n}^{(2)}}{\bi{3n}n}=\f1{10n+1}\prod\Sb j=1\\5\nmid j\endSb^{5n}\f{10n+j}j$$
is a $5$-adic integer. We are done. \qed

\proclaim{Lemma 3.3} Let $m\in\Z^+$ and $k,l,n\in\Z$. Then
$$\l\lfloor\f{kn+ln}m\r\rfloor-\l\lfloor\f{kn}m\r\rfloor-\l\lfloor\f{ln+1}m\r\rfloor
+\l\lfloor\f km\r\rfloor-\l\lfloor\f{k-1}m\r\rfloor\gs0.\tag3.2$$
\endproclaim
\Proof. If $m\nmid kn$, then
$$\l\lfloor\f{kn}m\r\rfloor+\l\lfloor\f{ln+1}m\r\rfloor=\l\lfloor\f{kn-1}m\r\rfloor+\l\lfloor\f{ln+1}m\r\rfloor
\ls\l\lfloor\f{(kn-1)+(ln+1)}m\r\rfloor.$$
If $m\nmid (ln+1)$, then
$$\l\lfloor\f{kn}m\r\rfloor+\l\lfloor\f{ln+1}m\r\rfloor=\l\lfloor\f{kn}m\r\rfloor+\l\lfloor\f{ln}m\r\rfloor
\ls\l\lfloor\f{kn+ln}m\r\rfloor.$$
Thus (3.2) holds when $kn$ or $ln+1$ is not divisible by $m$.

Now assume that $m\mid kn$ and $m\mid (ln+1)$. Clearly $m$ is relatively prime to $n$. Thus $m\mid k$ and hence
$$\align &\l\lfloor\f{kn+ln}m\r\rfloor-\l\lfloor\f{kn}m\r\rfloor-\l\lfloor\f{ln+1}m\r\rfloor
+\l\lfloor\f km\r\rfloor-\l\lfloor\f{k-1}m\r\rfloor
\\=&\f{kn}m+\f{ln+1}m-1-\f{kn}m-\f{ln+1}m+\f km-\l(\f km-1\r)=0.
\endalign$$

In view of the above, we have completed the proof of (3.2).
\qed

\medskip
\noindent{\it Proof of Theorem 1.4}. Clearly (1.12) holds if and only if
$(ln+1)\mid k\bi{kn+ln}{kn}$.

With the help of Lemma 3.3, for any prime $p$ we have
$$\align &\nu_p\(\f{k\bi{kn+ln}{kn}}{ln+1}\)=\nu_p\(\f{(kn+ln)!k!}{(kn)!(ln+1)!(k-1)!}\)
\\=&\sum_{j=1}^\infty\(\l\lfloor\f{kn+ln}{p^j}\r\rfloor
-\l\lfloor\f{kn}{p^j}\r\rfloor-\l\lfloor\f{ln+1}{p^j}\r\rfloor
+\l\lfloor\f k{p^j}\r\rfloor-\l\lfloor\f{k-1}{p^j}\r\rfloor\)\gs0.
\endalign$$
It follows that $ln+1$ divides $k\bi{kn+ln}{kn}$. We are done. \qed

\Ack. The author wishes to thank Dr. T. Amdeberhan, H. Q. Cao, I. Gessel and  H. Pan for helpful comments.

\enddocument

\bye